\documentclass[12pt,a5paper]{article}
\usepackage[a5paper,margin=6mm,tmargin=15mm]{geometry}

\title{Better buffers for patches in macroscale simulation of systems with microscale randomness}

\author{
J.~E. Bunder\thanks{School of Mathematical Sciences, University of Adelaide, South Australia~5005, Australia.  
\protect\url{mailto:judith.bunder@adelaide.edu.au}}
\and 
A.~J. Roberts\thanks{School of Mathematical Sciences, University of Adelaide, South Australia~5005, Australia.  
\protect\url{mailto:anthony.roberts@adelaide.edu.au}}
\and
I.~G. Kevrekidis\thanks{Department of Chemical and Biological Engineering and PACM,
Princeton University, Princeton, NJ~08544, USA.
\protect\url{mailto:yannis@Princeton.edu}}
}

\date{\today}

\usepackage{graphicx}
\usepackage{amsmath,url,doi}

\renewcommand{\vec}[1]{\text{\boldmath$#1$}}

\begin{document}

%\listofpdfcomments

\maketitle

\begin{abstract}
We consider one dimensional lattice diffusion model on a microscale grid with many discrete diffusivity values which repeat periodicially. 
Computer algebra explores how the dynamics of small coupled `patches' predict the slow emergent macroscale dynamics.
We optimise the geometry and coupling of patches by comparing the  macroscale predictions of the patch solutions with the macroscale solution on the infinite domain, which is derived for a general diffusivity period.
The results indicate that patch dynamics is a viable method for numerical macroscale modelling of microscale systems with fine scale roughness.
Moreover, the minimal error on the macroscale is generally obtained by coupling patches via `buffers' that are as large as half of each patch.  
\end{abstract}

\section{Introduction}
Across all fields of science and engineering there are examples of systems which are only realistically described by models on multiple spatial and temporal scales, where each scale is critical to the model and must be accounted for in any mathematical solution~\cite[e.g.]{dolbow2004, Horstemeyer2010, Dada2011}. 
Such multiscale systems are notoriously complex and computationally demanding, with numerical methods of solution generally requiring substantial reductions in complexity and careful control of errors. 
One significant aspect of multiscale modelling is how to transfer and represent information between different scales, particularly when the different scales are governed by different physics (for example, discrete and continuous, or deterministic and stochastic).
A typical scenario is that a macroscale model is required for practical applications, but essential to the system are fine details on one, or several, microscales which are substantially smaller than the macroscale of interest.
In developing the desired `coarse-grained' or macroscale model, major concerns are the efficiency of the the macroscale modelling and its accuracy compared to known, approximate macroscale closures. 
This article reports on an optimum design for one scheme---patch dynamics---for macroscale modelling.

Patch dynamics is an `equation-free' macroscale modelling technique in the sense that it makes no attempt to derive a macroscale closed form from the original microscale model~\cite[for reviews]{Givon04, Hyman05, Li2007, Kevrekidis09, Samaey10}.  
Numerical solutions are obtained on demand by solving the microscale model on small discrete patches, with each patch centred about one macroscale grid point.
A major advantage of this technique is that it is the original microscale model which is solved numerically, and not an approximation.
Unlike many other macroscale modelling techniques~\cite[e.g.]{Pavliotis2008, Degenhard2006}, no assumptions are made regarding the relative importance of certain components of the microscale model since patch dynamics requires no such simplification. 
The main complications of patch dynamics are (a) how to couple the discrete patches across space so that the unevaluated regions between the patches are best accounted for, and (b) how to determine an appropriate placement and size of patch which captures sufficient microscale detail. 

A full implementation of patch dynamics involes macroscale modelling in both space and time~\cite{Runborg2002, Moller2005, Givon04, Samaey05}. 
Here we concentrate on the spatial aspects of the modelling.   
Spatial patch dynamics, also known as the gap-tooth scheme, was successfully applied to various models, including Burgers' equation~\cite{Roberts01}, a generalized advection-diffusion equation~\cite{Roberts03} and a Ginzburg--Landau model~\cite{Roberts11}, and the macroscale modelling was shown to be largely independent of both patch width and patch coupling conditions~\cite{Roberts07}. 
However, these models all describe systems where the microscale structures vary smoothly.
Here we explore an infinite one dimensional lattice diffusion model, specified in Section~\ref{sec:model}, but in contrast to earlier studies we invoke microscale detail characterised by non-smooth discrete diffusivity with rough spatial variability.

The importance of microscale diffusion on a lattice is that it is in the same universality class as a host of other microscale models.
We thus contend that our results apply to a wide variety of microscale systems whose emergent macroscale dynamics are that of diffusion.
We also expect the results suggest first approximations to suitable patch construction and patch coupling conditions for other more complicated dynamical systems. 

To ascertain the ability of patch dynamics to model fine scale spatial roughness, our diffusion model contains $K\geq 2$ independent discrete diffusivity values which repeat periodicially on an infinite microscale lattice.
We previously considered a similar model~\cite{Bunder2012}, although only for the case of two alternating diffusivities, $K=2$\,.
Analogously Knapek~\cite{Knapek98} developed a homogenization-based coarse-graining method for a two dimensional diffusion model with variable diffusion which showed varying degrees of success. 
We generalise previous research of two discrete diffusivity values to many discrete diffusivity values and compare the patch dynamics macroscale modelling with the macroscale closure across the complete domain, derived in Section~\ref{sec:infinite}. 
For an accurate macroscale closure we must ensure that the symmetries of the original microscale model, defined by the fine microscale structure, are preserved by the macroscale modelling.
The microscale model has translational symmetry since any shift of the $K$ diffusivity values along the infinite microscale lattice, while maintaining their original order, does not affect the microscale solution. 
In addition, the microscale model has reflection symmetry since reversing the order of the $K$ diffusivity values has no effect on the  microscale solution.
Section~\ref{sec:patch} constructs the patch dynamics method for our diffusion model and discusses the need for an ensemble average to avoid asymmetric effects~\cite{Moller2005}.
Section~\ref{sec:optimise} compares the emergent patch dynamics closure with the complete domain solution. 
By minimising the errors, we determine the optimal patch coupling.

Section~\ref{sec:exact} discusses the special cases where appropriately chosen patch geometry and coupling between patches fully preserve the symmetries associated with the fine microscale structure without requiring an ensemble average.
For these cases the desired macroscale closure is fully reproduced, to any desired order of accuracy. 
This shows the importance of careful patch construction for models with fine scale detail when using patch dynamics for macroscale modelling.

\section{Discrete microscale diffusion equation}
\label{sec:model}

We consider discretised diffusion on an infinite one dimensional microscale lattice with lattice spacing~$h$,
\begin{equation}
\dot{u}_i(t)=\kappa_{i}\left[u_{i+1}(t)-u_i(t)\right]/h^2+\kappa_{i-1}\left[u_{i-1}(t)-u_i(t)\right]/h^2\,,\label{eq:de}
\end{equation}
where $\dot{u}_i(t)$ is the time derivative of the field~$u_i(t)$ at lattice point~$i$ and
 $\kappa_{i}=\kappa_{i\operatorname{mod}K}$ is the diffusivity value at half lattice point $i+1/2$\,. 
There are $K\geq 2$ independent diffusivity values $(\kappa_1,\kappa_2,\ldots,\kappa_K)$ which repeat periodically on the microscale lattice with period $K$.
The ultimate aim is to simulate the dynamics of equation~\eqref{eq:de}, not for all~$i$, but on a significantly larger macroscale grid with spacing $H\gg h$\,. 
We define the macroscale lattice points on the macroscale grid as~$X_j$ for integer~$j$ and seek the dynamics of a discrete macroscale field~$U_j(t)$ for all~$j$.
The macroscale field~$U_j(t)$ is to describe the large scale dynamics of the microscale field~$u_i(t)$ defined by equation~\eqref{eq:de}.

\begin{figure}
\centering\includegraphics{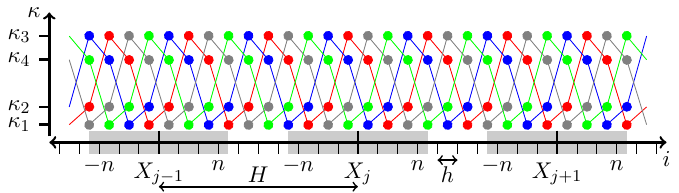}
\caption{The microscale lattice is indexed by~$i$ and indicated by the fine short ticks on the horizontal axis with spacing~$h$.
The macroscale lattice is indexed by~$X_j$ and indicated by the thick long ticks on the horizontal axis  with spacing~$H$.
We construct the $j$th~patch of width~$(2n+1)h$ about the macroscale lattice point~$X_j$ for all~$j$, indicated by the shaded rectangles.
Here we represent a microscale model with $K=4$ given diffusivity values~$\kappa_{1,2,3,4}$ and we choose patch half-width $n=3$\,. 
The ensemble contains $2K=8$ configurations; the four with translation symmetry are illustrated, starting from the leftmost marker: $(\kappa_1,\kappa_2,\kappa_3,\kappa_4)$ (grey), $(\kappa_2,\kappa_3,\kappa_4,\kappa_1)$ (red), $(\kappa_3,\kappa_4,\kappa_1,\kappa_2)$ (blue), $(\kappa_4,\kappa_1,\kappa_2,\kappa_3)$ (green).
The remaining four configurations which complete the ensemble are reflections of the illustrated configurations;  $(\kappa_4,\kappa_3,\kappa_2,\kappa_1)$\,, $(\kappa_1,\kappa_4,\kappa_3,\kappa_2)$\,, $(\kappa_2,\kappa_1,\kappa_4,\kappa_3)$\,,  and $(\kappa_3,\kappa_2,\kappa_1,\kappa_4)$\,. }
\label{fig:lattice}
\end{figure}

As illustrated in Figure~\ref{fig:lattice}, we construct small patches of width $(2n+1)h<H$ for positive integer patch half-width~$n$, centered about each macroscale lattice point~$X_j$. 
We use the microscale equation~\eqref{eq:de} for~$u_i$, but only within each patch rather than across the entire microscale domain.
To distinguish the fields on the complete microscale domain~$u_i$ from the patch dynamics fields in the $j$th~patch we represent the latter by~$u_{j,i}$ for $i=0,\pm1,\ldots,\pm n$\,.
To obtain a well-posed problem for the patch dynamics we must supply appropriate boundary conditions for each patch.
These boundary conditions are generally referred to as coupling conditions since their purpose is to couple adjacent patches by extrapolating across the unevaluated space between patches, thus providing a solution which approximates that of the complete domain. 
The macroscale field~$U_j$ is obtained by some averaging over the $u_{j,i}$~field within the $j$th~patch.  
Section~\ref{sec:patch} presents details concerning the patch construction, coupling conditions and obtaining the macroscale solution.

The aim is to show that, by providing suitable patch coupling conditions and an appropriate averaging over the microscale fields within each patch, we are able to obtain a description of the large scale dynamics of the system defined by equation~\eqref{eq:de} while substantially reducing the numerical cost of evaluating the solution over the entire microscale domain.
Section~\ref{sec:infinite} derives the long scale evolution of equation~\eqref{eq:de} across the infinite microscale domain, 
and Section~\ref{sec:optimise} shows that patch dynamics is able to capture this emergent behaviour. 
 Section~\ref{sec:optimise} also discusses how the success of patch dynamics is dependent on choosing appropriate patch geometry and coupling conditions for a given microscale structure. 

The fine microscale structure of the model, defined by the roughly varying diffusivities~$\kappa_i$, 
complicates the problem for several reasons. 
Fine scale roughness across the microscale lattice produces a $u_i$~field solution which is similarly rough, possibly compromising the macroscale modelling of the emergent dynamics.
Care must be taken so that important details concerning this microscale structure are not lost. 
For example, not all diffusivity values will be sampled if the patch is smaller than the period of the diffusivities, that is $(2n+1)<K$, and the patches may destroy symmetries inherent in the original model.
Section~\ref{sec:optimise} shows that patch dynamics is able to overcome these issues and fully describe the emergent dynamics to desired accuracy, provided a suitable patch geometry is chosen with complementary patch coupling conditions.
The issue of destroyed symmetries arises because, in the original model there is translational and reflection symmetry across the complete microscale, but within a patch of arbitrary but finite size there are generally no such symmetries.
While specific patch geometries do eliminate this problem, there are circumstances where we may not be able to choose the optimal geometry.
Therefore, to reintroduce lost symmetries into the patches we evaluate over an ensemble of diffusivity values, rather than one specific configuration~\cite{Bunder2012}. 
The full ensemble contains~$2K$ configurations comprising $K$~translations of the diffusivity values and $K$~reflections of the diffusivity values.
Section~\ref{sec:infinite} shows that the emergent dynamics is dependent on the order of the $K$ diffusivities so we do not include all possible diffusivity permutations in the full ensemble but only those which preserve their order. 
For $K\geq 3$ the $2K$ configurations are unique, 
but for $K=2$ there are only two unique configurations.
For example, Figure~\ref{fig:lattice} shows the four translated configurations for $K=4$.
By ensemble averaging, patch dynamics is able to adequately capture the large scale emergent dynamics of the original problem~\eqref{eq:de}, even when the patch geometry is not optimal.

\section{Emergent evolution derived on the complete microscale}
\label{sec:infinite}

We derive the emergent large scale dynamics of the microscale discrete diffusion equation~\eqref{eq:de} for~$K$ given  diffusivities~$\kappa_i$ on the complete microscale domain. 
The resultant description of the dynamics provides a basis for measuring errors in the subsequent patch modelling.

In a dissipative system such as~\eqref{eq:de} the long term, emergent dynamics are characterised by the smallest magnitude eigenvalue~$\lambda_0$, that is, $u_i\sim e^{\lambda_0t/h^2}$ after a sufficiently long time~$t$.  
The emergent evolution of the macroscale field~$U_j$ must therefore be correspondingly $\dot{U}_j=\lambda_0U_j/h^2$. 
Here we construct an analytic expression for the complete domain emergent evolution, $\dot{U}_j=g_0(\vec U)$, which Section~\ref{sec:optimise} compares with the evolution obtained from patch dynamics in order to establish the accuracy of the patch dynamics numerical scheme and to determine optimal parameters.

We partition the domain into elements covering one period of the diffusivity, and then consider the dynamics within each cell but coupled to neighbouring cells.
Equation~(\ref{eq:de}) in matrix form is then written as $\dot{{\vec u}}=M{\vec u}/h^2$, for field vector ${\vec u}=(u_i,u_{i+1},\ldots,u_{i-1+K})$ and where the only nonzero elements of the $K\times K$ operator matrix~$M$ are
\begin{align}
M_{i,i}&=-\kappa_{i-1}-\kappa_i\,,  &M_{i,i+1}&=M_{i+1,i}=\kappa_i\text{ for } i<K,\nonumber\\
M_{1,K}&=\kappa_K\varepsilon ^{-K}, &M_{K,1}&=\kappa_K\varepsilon ^{K}.
\end{align}
In these expressions, the operator~$\varepsilon$ is a microscale step operator which shifts a field by one microscale lattice spacing~$h$; that is, $\varepsilon u_i=u_{i+1}$ and its inverse $\varepsilon^{-1} u_i=u_{i-1}$\,. 
The microscale step operator~$\varepsilon$ in the matrix~$M$ allow us to reduce the original problem of solving a potentially infinite number of variables~$u_i$ to a problem with just~$K$ variables.
The cost is that the shift operator occurs in the description of the \(K\)~variables so that we implicitly cater for an infinite domain.

The dynamics of the matrix form of the microscale system is characterised by the $K$~eigenvalues~$\lambda$ of the $K\times K$~matrix~$M$ which satisfy the $K$th~order characteristic equation 
\begin{equation}
c_0({\vec \kappa})+c_1({\vec \kappa})\lambda+c_2({\vec \kappa})\lambda^2+\cdots+c_{K}({\vec \kappa})\lambda^{K}=0\,,\label{eq:chareq}
\end{equation}
for coefficients $c_q$ which are functions of the diffusivities ${\vec \kappa}=(\kappa_1,\ldots,\kappa_K)$. 
The coefficients in the characteristic equation~\eqref{eq:chareq} are 
\begin{align}
c_0({\vec \kappa})={}&-\left(\varepsilon^{K/2}-\varepsilon^{-K/2}\right)^2\kappa_g^K\,,\label{eq:cs}\\
c_q({\vec \kappa})={}&\frac{\kappa_g^K}{q}\sum_{m_1=1}^K\sum_{m_2=1}^{(K-q+1)}\sum_{m_3=1}^{(K-q+2-m_2)}
\times\cdots \nonumber\\
{}&\cdots\times\sum_{m_q=1}^{(K-1-m_2-\cdots-m_{q-1})}\frac{m_2\cdots m_q(K-m_2-\cdots -m_q)}
{\kappa_{m_1}\kappa_{m_1+m_2}\cdots\kappa_{m_1+m_2+\cdots+m_q}}\,,\nonumber
\end{align}
for $1\leq q\leq K$\,. 
Two basic cases are $c_K=1$  and $c_1=K^2\kappa^K_g/\kappa$ where
\begin{equation}
\kappa_g=\left(\prod_{i=1}^K\kappa_i\right)^{1/K}
\quad\text{and}\quad 
\kappa=K\left(\sum_{i=1}^K\frac{1}{\kappa_i}\right)^{-1}
\end{equation}
are the geometric and harmonic means of the diffusivities, respectively.

We are only interested in the smallest magnitude eigenvalue $\lambda_0$ which describes the emergent dynamics of the system, rather than all $K$~eigenvalues.
In a diffusion model we expect $|\lambda_0|\ll 1$ and therefore
 we approximate the characteristic equation~\eqref{eq:chareq} by its quadratic truncation,
\begin{equation}
c_0({\vec \kappa})+c_1({\vec \kappa})\lambda_0+c_2({\vec \kappa})\lambda_0^2=0\,.
\end{equation} 
% $c_0=-2K\kappa_g^K\sum_{p=1}^K\frac{(K+p-1)!}{(K-p)!(2p)!}\delta^{2p}$
The quadratic coefficient is
\begin{equation}
c_2({\vec \kappa})=\frac{\kappa_g^K}{2}\sum_{m_1=1}^K\sum_{m_2=1}^{K-1}\frac{m_2(K-m_2)}
{\kappa_{m_1}\kappa_{m_1+m_2}}\,.
\end{equation}
Unlike the linear coefficient~$c_1$ and the constant coefficient~$c_0$, the quadratic coefficient is generally dependent on the order of the diffusivities on the microscale lattice.
For example, when $K=4$
\begin{equation}
c_2({\vec \kappa})=3\sum_{i=1}^4\kappa_i\kappa_{i+1}+4\sum_{i=1}^2\kappa_i\kappa_{i+2}
\end{equation}
is different for configurations of the diffusivities which are not translations or reflections of each other, 
that is,
${\vec \kappa}=(\kappa_1,\kappa_2,\kappa_3,\kappa_4)$ and ${\vec \kappa}=(\kappa_1,\kappa_3,\kappa_2,\kappa_4)$ do not have the same coefficient~$c_2({\vec \kappa})$\,, whereas the three configurations ${\vec \kappa}=(\kappa_1,\kappa_2,\kappa_3,\kappa_4)$,  ${\vec \kappa}=(\kappa_2,\kappa_3,\kappa_4,\kappa_1)$ and ${\vec \kappa}=(\kappa_1,\kappa_4,\kappa_3,\kappa_2)$ have the same coefficient~\(c_2({\vec \kappa})\)\,.
For this reason, when we take the ensemble average with patch dynamics we only consider those configurations which preserve the order of the diffusivities; that is, only translations and reflections of the original order of diffusivity values as defined by the original microscale model~\eqref{eq:de}. 
This dependence on the order of the diffusivities only becomes important when period $K\geq 4$ since for $K=2$ and~$3$ all configurations are reflections or translations of each other.
The dependence of~$c_2$ on the ordering of diffusivities for $K\geq 4$ shows how the detailed microscale structure of system affects the solution on the macroscale and implies that a macroscale model which is only dependent on mean diffusivities may not be sufficiently accurate---depending upon the accuracy required.
We show in Section~\ref{sec:optimise} that patch dynamics is able to account for different configurations of the microscale structure.   

For a quadratic approximation of the characteristic equation, the smallest magnitude solution~$\lambda_0$ produces the macroscale large scale emergent evolution
\begin{align}
\dot U_j={}&g_0(\vec U)=\lambda_0U_j/h^2
\nonumber\\
={}&-\frac{K^2\kappa_g^K}{2h^2\kappa c_2({\vec \kappa})}\left[1-\sqrt{1+
4\kappa c_2({\vec \kappa})\left(\varepsilon^{K/2}-\varepsilon^{-K/2}\right)^2/K^2}\right]U_j\nonumber\\
={}&\frac{\kappa}{h^2}\left[\delta^2+\left(\frac{K^2-1}{12}-\frac{c_2({\vec \kappa})\kappa^2}{K^2\kappa_g^K}\right)
\delta^4\right]U_j+\mathcal O(\delta^6)\,.\label{eq:lambda0}
\end{align}
The $O(\delta^6)$ error in~$g_0$ is due to our quadratic approximation of the characteristic equation, with the exception of $K=2$ where the quadratic truncation is exact. 
The expansion of the square root uses
\begin{equation}
\left(\varepsilon^{K/2}-\varepsilon^{-K/2}\right)^2
=K^2\delta^2\left[1+(K^2-1)/12\delta^2\right]
+\mathcal O(\delta^6),
\end{equation}
where, for $K>2$, we should only expand the square root up to $O(\delta^6)$ error since this corresponds to the error from the quadratic truncation.
The microscale difference operator $\delta=\varepsilon^{1/2}-\varepsilon^{-1/2}$ is treated as small in the sense that $\delta^2 U_j=(\varepsilon+\varepsilon^{-1}-2)U_j\ll U_j$ since $\delta^2U_j$ measures a microscale variation of the slow macroscale solution~$U_j$, and this variation must be small. 
This supports our previous assumption that $|\lambda_0|\ll 1$\,.
Section~\ref{sec:optimise} uses equation~\eqref{eq:lambda0} to determine the accuracy of the macroscale evolution obtained from patch dynamics.

\section{Patch dynamics}
\label{sec:patch}

Section~\ref{sec:construct} provides details for constructing discrete patches on the complete microscale domain and highlights issues particular to systems which contain fine microscale detail,  such as equation~\eqref{eq:de}. 
We then discuss the method for obtaining the macroscale field solutions and define the coupling conditions. 
%However, one must be aware that there are no strict rules on how to define coupling conditions or extract the macroscale solution.
%A good coupling condition is simply any condition which is easy to apply numerically and couples the patches in some reasonable way so that there is good extrapolation across the unevaluated space between the patches. 
%Similarly, the macroscale field solution~$U_j$ at $X_j$ should be projected from the microscale field solutions~$u_{j,i}$ within the $j$th~patch in some numerically simple manner, but there is no formal theorem which tells us precisely how this should be done. 
Our method for obtaining the macroscale fields in Section~\ref{sec:extract} and  the coupling conditions in Section~\ref{sec:coupling} describe one possible patch dynamics macroscale modelling method amongst many.
The test of a good patch dynamics method is simply how well it describes known emergent dynamics. 
In Section~\ref{sec:optimise} we show that our method does effectively reproduce the known emergent dynamics of the original microscale system~\eqref{eq:de}, 
as derived in Section~\ref{sec:infinite}.

\subsection{Patch construction}
\label{sec:construct}
We construct the $j$th~patch of width~$(2n+1)h$, for integer patch half-width~$n$ and microscale lattice spacing~$h$, about the macroscale lattice point~$X_j$, for all~$j$, as indicated in grey in Figure~\ref{fig:lattice}. 
This division of the complete microscale domain into discrete patches may introduce some unwanted features into the macroscale modelling which must be removed in order to obtain the best possible macroscale closure.
For example, if the patches are positioned such that diffusivity~$\kappa_1$ always appears on the left edge of each patch, this gives~$\kappa_1$ a significance which does not exist in the original problem.
Also, within one patch some diffusivity values may appear more often than others, which again gives greater significance to some microscale structures over others.
These issues, as briefly discussed in Section~\ref{sec:model}, are due to the patches destroying symmetries of the original complete domain model~\eqref{eq:de}. 
If these lost symmetries are not taken into account when constructing the patch dynamics method, then undesirable terms tend to arise in the macroscale modelling.
For example, drift terms, which describe differences in the left and right moving flows, are expected in an advection-diffusion model, but not  in a diffusion model, and yet they do appear in the patch dynamics emergent evolution of a diffusion model if symmetries are not correctly accounted for.

To reinstate the correct symmetry our patch dynamics method generally requires an average over an ensemble of diffusivity configurations~\cite{Moller2005}.
This ensures that all diffusivities are treated equally and the model is independent of how the patches are positioned on the complete domain.  
The original microscale model~\eqref{eq:de} is symmetric under translations and reflections of the diffusivity values, but does not permit changing the order of the diffusivity values.
Section~\ref{sec:infinite} derives the macroscale evolution on the complete domain and shows that it also has translation and reflection symmetry. 
Therefore, we construct an ensemble of diffusivity configurations where each configuration describes either a translation or reflection of the original diffusivity configuration ${\vec \kappa}=(\kappa_{1},\ldots,\kappa_{K})$\,, while maintaining the same ordering of diffusivity values.  
For $K$ diffusivity values the ensemble contains a total of~$2K$~configurations consisting of $K$ translations and $K$ reflections of the original configuration~${\vec \kappa}$\,.
We represent one configuration by the subscript~$e$,  ${\vec \kappa}_e=(\kappa_{1,e},\ldots,\kappa_{K,e})$, and the discrete diffusion  equation~\eqref{eq:de} on the $j$th~patch with configuration~$e$ is 
 \begin{equation}
\dot{u}_{j,i,e}(t)=\kappa_{i,e}\left[u_{j,i+1,e}(t)-u_{j,i,e}(t)\right]/h^2
+\kappa_{i-1,e}\left[u_{j,i-1,e}(t)-u_{j,i,e}(t)\right]/h^2\,,\label{eq:depatch}
\end{equation}
for $i=0,\pm 1,\ldots,\pm (n-1)$ and $e=1,2,\ldots,2K$\,. 
Our patch dynamics method is to solve equation~\eqref{eq:depatch} for all fields~$u_{j,i,e}$ in all patches~$j$ and for all ensembles~$e$, with patch coupling conditions defined in Section~\ref{sec:coupling}. 
The macroscale field solution~$U_j$ at each~$X_j$ is defined in Section~\ref{sec:extract}.
Section~\ref{sec:procedure} provides a complete description of the patch dynamics modelling procedure, including step-by-step instructions for applying the coupling conditions and obtaining the time dependent macroscale field solutions~$U_j(t)$.

\subsection{Extracting a macroscale field}
\label{sec:extract}

\begin{figure}
\centering\includegraphics{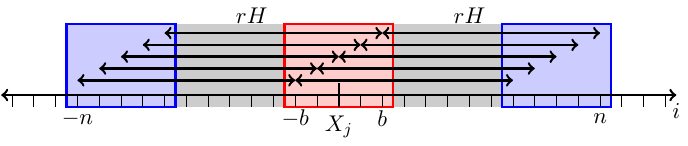}
\caption{A closeup view of one patch, similar to those shown in Figure~\ref{fig:lattice} but with patch half-width $n=12$ and buffer half-width $b=2$\,. 
The outlined shaded region in the centre of the patch containing $i=0,\pm1,\ldots,\pm b$ defines the core and is required for the amplitude condition~\eqref{eq:amp}. 
The two outlined shaded regions on the ends of the patch containing $i=\pm(n-2b),\ldots,\pm n$ are the buffers and are required for the coupling conditions~\eqref{eq:cc}. 
Each point in each buffer is a distance of $rH=(n-b)h$ from one point in the core.}
\label{fig:buffers}
\end{figure}

Once equation~\eqref{eq:depatch} within the $j$th~patch is solved for~$u_{j,i,e}$ with appropriate coupling conditions, we need to extract a macroscale solution field from these~$u_{j,i,e}$. 
We propose a macroscale field solution, generally called the amplitude condition, which is both an ensemble average and an average over a number of microscale fields $u_{i,j,e}$ in the centre, or `core', of the $j$th patch,
\begin{equation}
U_j(t)=U(X_j,t)=\left\langle\sum_{i=-b}^b\frac{u_{j,i,e}}{2b+1}\right\rangle,\label{eq:amp}
\end{equation}
where the angle brackets represent the ensemble average over all configurations~$e$, and integer~\(b\), which we call the buffer half-width, satisfies $0\leq b<n$\,.
The amplitude condition requires an ensemble average because we only want one macroscale field solution at each~$X_j$, not a macroscale solution for each configuration~$e$. Figure~\ref{fig:buffers} illustrates one patch and highlights the core region.

Some patch dynamics methods do not suggest an average over the patch core and extract a macroscale solution by simply using the microscale solution at the centre of the patch~\cite[e.g.]{Roberts01, Roberts11}, corresponding to $b=0$ in equation~\eqref{eq:amp}. 
However, these methods are typically for systems with smooth microscale dynamics, resulting in a smooth microscale solution which varies only slightly within the patch.
Given our rough microscale structure, which produces a rough microscale field, an average over the microscale solutions within the patch core seems appropriate as it should reduce wild fluctuations which do not adequately reflect the larger macroscale dynamics. 
Section~\ref{sec:optimise} considers how the choice of buffer half-width~$b$ affects the large scale emergence of the macroscale dynamics, and determines how to best choose~$b$ based on the underlying microscale structure of the original model~\eqref{eq:de}.

\subsection{Coupling conditions}
\label{sec:coupling}

To solve equation~\eqref{eq:depatch} for~$u_{j,i,e}$ within the $j$th~patch with $i=0,\pm1,\ldots,\pm n$ for one particular configuration~$e$, we need to define two coupling conditions.
The coupling conditions for the $j$th~patch are derived from an interpolation of the macroscale field~$U_j$, as defined by the amplitude condition~\eqref{eq:amp}, and its nearest macroscale neighbours, $U_{j\pm 1}$, $U_{j\pm 2}$,\ldots\,, and are designed to constrain the microscale fields~$u_{j,i,e}$ on or near the edges of each patch.
Many patch dynamics coupling conditions only apply to the patch end points $i=\pm n$~\cite[e.g.]{Roberts01, Roberts11}. 
However, for roughly structured microscales, fixing the microscale fields on the patch end points,~$u_{j,\pm n,e}$, is not practical.
The rough structure of the microscale produces solutions of~$u_{j,i,e}$ with a similarly rough structure which tends to be compromised by strict requirements at a single point such as $i=\pm n$.
Therefore, like the amplitude condition~\eqref{eq:amp}, we use an average over several\slash many microscale grid points to define a coupling condition.
These averaged regions, the so-called `buffers', avoid strictly fixing the microscale field at the patch end points $i=\pm n$ and provide  some limited `freedom' on the patch edges. 
Samaey et al.~\cite{Samaey2004} showed that it is possible to define quite arbitrary patch coupling conditions, provided suitably large buffers are chosen to shield the patch core from unwanted consequences of these coupling conditions over short evolution times. 
The buffers in our model serve quite a different purpose than the buffers of Samaey et al.~\cite{Samaey2004}, but in both cases, when chosen appropriately, offer significant advantages, as shown in Section~\ref{sec:optimise}.

We derive coupling conditions compatible with the amplitude condition~\eqref{eq:amp} as then some beautiful results follow in certain circumstances.
Our coupling conditions invoke the microscale and macroscale step operators, 
$\varepsilon$~and~$\bar{\varepsilon}$, respectively.
Section~\ref{sec:infinite} defined the microscale step operator~$\varepsilon$ to shift the microscale field by one microscale lattice spacing~$h$.
Here we have the added complication of patches and multiple configurations, but the shift~$\varepsilon$ still operates in the same way; that is 
$\varepsilon^{\pm 1} u_{j,i,e}=u_{j,i\pm 1,e}$\,. 
To shift the right-hand side of~\eqref{eq:amp} to the right patch buffer we apply~$\varepsilon^{n-b}$ to both side of the equation, whereas to shift to the left patch buffer we apply the inverse~$\varepsilon^{-(n-b)}$:
\begin{equation}
\varepsilon^{\pm(n-b)}U_{j}=\left\langle\sum_{i=n-2b}^{n}\frac{u_{j,\pm i,e}}{2b+1}\right\rangle.\label{eq:temp}
\end{equation}
To write the left-hand side of the above equation as something meaningful we define~$r$ such that $rH=(n-b)h$\,,
which equates the $(n-b)$~microscale steps to a fractional macroscale step of~$r$.
We also define the macroscale step operator~$\bar{\varepsilon}$ which shifts the macroscale field by one macroscale lattice spacing~$H$, $\bar{\varepsilon}^{\pm}U_{j,e}=U_{j\pm 1,e}$. 
Thus, in equation~\eqref{eq:temp}, 
 $\varepsilon^{\pm (n-b)}U_{j}=\bar{\varepsilon}^{\pm r}U_{j}$\,.
Finally, we introduce the macroscale difference and mean operators, $\bar{\delta}= \bar{\varepsilon}^{1/2}-\bar{\varepsilon}^{-1/2}$ and
$\bar{\mu}=(\bar{\varepsilon}^{1/2}+\bar{\varepsilon}^{-1/2})/2$, respectively,
rewrite~$\bar{\varepsilon}^{\pm r } U_{j}$ in terms of $\bar{\mu}$ and $\bar{\delta}$, and expand in powers of $\bar{\delta}^2$, 
which are small when operating on~$U_{j}$ since amplitude~$U_{j}$ varies slowly in space.
Thus we obtain the coupling conditions defined in terms of averages over buffers:
\begin{align}
\sum_{i=n-2b}^n \frac{u_{j,\pm i,e}}{2b+1}
={}&U_{j}+\sum_{k=1}^{\Gamma}\left(\prod_{l=0}^{k-1}(r^2-l^2)\right)\gamma^k\frac{\pm (2k/r)\bar{\mu}\bar{\delta}^{2k-1}+\bar{\delta}^{2k}}{(2k)!}U_{j}
\nonumber\\{}&
+\mathcal{O}(\gamma^{\Gamma+1}),\label{eq:cc}
\end{align}
where $\Gamma$ is some cutoff for the expansion in terms of the coupling strength parameter~$\gamma$.
In the coupling conditions we do not have an ensemble average because we need one coupling condition at each patch edge for each configuration $e$.
Figure~\ref{fig:buffers} illustrates the shifts of $rH=(n-b)h$ from the patch core to the patch buffers which describe the operation of 
$\bar{\varepsilon}^{\pm r}=\varepsilon^{\pm(n-b)}$ on equation~\eqref{eq:amp} to obtain~\eqref{eq:temp} and ultimately the coupling conditions~\eqref{eq:cc}.

The coupling strength~$\gamma$ is an artificial parameter which controls the coupling between patches. 
It is important for establishing theoretical support for patch dynamics.
Detailed discussions on the theoretical support for patch dynamics through the introduction of the coupling strength~$\gamma$ are provided elsewhere~\cite{Roberts03, Roberts07, Samaey10}. 
We greatly summarise the main points here.
For the unphysical case of no coupling $\gamma=0$\,, since all eigenvalues of the microscale system on the $j$th~patch have negative real-parts except for one zero eigenvalue, a slow manifold exists of dimensionality equal to the number of patches. 
This slow manifold persists to non-zero coupling~\(\gamma\).
The whole system dynamics is attracted to this slow manifold on a cross-patch diffusion time, and thereafter the system evolves slowly to form the macroscale dynamics of interest.
We approximate the slow manifold by a power series in coupling~\(\gamma\).
Although the radius of convergence is unknown, we conjecture that evaluation at the physical case of full coupling $\gamma=1$ is sufficiently accurate to form useful models as it does in other examples~\cite[e.g.]{Roberts03, Roberts07}.
The error in the coupling conditions of order~$\gamma^{\Gamma+1}$ is equivalent to errors of order~$\bar{\delta}^{2(\Gamma+1)}U_j$ and~$\bar{\mu}\bar{\delta}^{2\Gamma+1}U_j$, which are small provided the macroscale field~$U_j$ varies slowly across the macroscale spacing~$H$.

The coupling conditions are defined in terms of the macroscale mean and difference operators~$\bar{\mu}$ and~$\bar{\delta}$, but the complete domain evolution~\eqref{eq:lambda0} is dependent on the microscale operator~$\delta$.
To effectively compare the evolution obtained from patch dynamics and the evolution on the complete domain we must derive the relationship between the microscale (no bar) and macroscale (bar) operators.  
To do this we also need to define a microscale mean operator $\mu=(\varepsilon^{1/2}+\varepsilon^{-1/2})/2$. 
To convert a macroscale operator into microscale operators we write the macroscale operator in terms of the macroscale step operator $\bar{\varepsilon}$.
Then, we substitute the previously derived relationship between the two step operators, $\bar{\varepsilon}^{\pm 1}=\varepsilon^{\pm(n-b)/r}$\,, to obtain an equation in terms of the microscale operator $\varepsilon$. 
For our purposes it is more convenient to use mean and difference operators rather than step operators and so we then substitute $\varepsilon^{\pm1}=1\pm\mu\delta+\delta^2/2$ to obtain the original macroscale operator in terms of microscale mean and difference operators. 

For example, we expand the macroscale operator~$\bar{\delta}^2$ in terms of~$\bar{\varepsilon}$, $\bar{\delta}^2=\bar{\varepsilon}+\bar{\varepsilon}^{-1}-2$ and then substitute $\bar{\varepsilon}^{\pm 1}=\varepsilon^{\pm(n-b)/r}$\,. 
This gives the macroscale operator $\bar{\delta}^2$ in terms of the microscale operator $\varepsilon$. 
Then we substitute $\varepsilon^{\pm1}=1\pm\mu\delta+\delta^2/2$\,, and on performing a binomial expansion
\begin{equation}
\bar{\delta}^2=\sum_{l=1}^{\infty}\frac{1}{l!}\left[\prod_{k=0}^{l-1}\left(\frac{n-b}{r}-k\right)\right]\left[(\mu\delta+\delta^2/2)^l+(-\mu\delta+\delta^2/2)^l\right].\label{eq:convert1}
\end{equation}
One similarly derives the inverse relationship, from microscale operator~$\delta^2$ to macroscale operators $\bar{\delta}$~and~$\bar{\mu}$,
\begin{equation}
\delta^2=\sum_{l=1}^{\infty}\frac{1}{l!}\left[\prod_{k=0}^{l-1}\left(\frac{r}{n-b}-k\right)\right]\left[(\bar{\mu}\bar{\delta}+\bar{\delta}^2/2)^l+(-\bar{\mu}\bar{\delta}+\bar{\delta}^2/2)^l\right].
\label{eq:delbar2del}
\end{equation}
The conversion of $\bar{\mu}\bar{\delta}$ into microscale operators is also useful. 
Following a similar procedure as above, but with $\bar{\mu}\bar{\delta}=(\bar{\varepsilon}-\bar{\varepsilon}^{-1})/2$\,,
\begin{align}
\bar{\mu}\bar{\delta}&=\sum_{l=1}^{\infty}\frac{1}{2l!}\left[\prod_{k=0}^{l-1}\left(\frac{n-b}{r}-k\right)\right]\left[(\mu\delta+\delta^2/2)^l-(-\mu\delta+\delta^2/2)^l\right],\nonumber\\
\mu\delta&=\sum_{l=1}^{\infty}\frac{1}{2l!}\left[\prod_{k=0}^{l-1}\left(\frac{r}{n-b}-k\right)\right]\left[(\bar{\mu}\bar{\delta}+\bar{\delta}^2/2)^l-(-\bar{\mu}\bar{\delta}+\bar{\delta}^2/2)^l\right].\label{eq:convert2}
\end{align}

The complete domain evolution~\eqref{eq:lambda0} derived in Section~\ref{sec:infinite} is correct to $\mathcal{O}(\delta^4)$.  
To obtain the same order of accuracy from our patch dynamics modelling we choose cutoff~$\Gamma=2$ in the coupling conditions~\eqref{eq:cc} and when converting between macroscale and microscale operators we require up to $l=4$ in the infinite sums in equations~\eqref{eq:convert1}--\eqref{eq:convert2}. 

\subsection{Patch dynamics macroscale modelling procedure}
\label{sec:procedure}
Here we describe, step-by-step, how one obtains a macroscale closure of  the microscale model~\eqref{eq:de} using  patch dynamics macroscale modelling, as defined by equation~\eqref{eq:depatch}, the coupling conditions~~\eqref{eq:cc} and the amplitude condition~\eqref{eq:amp}.
Say we wish to numerically simulate equation~\eqref{eq:de} with initial conditions on the complete domain of $u_{i}(0)=u^0_i$ for all~$i$ and periodic diffusivities ${\vec \kappa}=(\kappa_1,\ldots,\kappa_K)$ at macroscale lattice points~$X_j$.  
We firstly determine all $2K$ configurations~${\vec \kappa}_e$, construct patches about each~$X_j$, and redefine the microscale initial conditions as initial conditions within each patch~$u_{j,i,e}(0)$. 
Then we use the amplitude condition~\eqref{eq:amp} and the patch microscale initial conditions~$u_{j,i,e}(0)$ to determine the macroscale initial conditions~$U_{j}(0)$ for all~$j$. 
If we have evaluated the macroscale solutions $U_{j}(t_{m-1})$ at time step~$t_{m-1}$ for all~$j$, then the patch dynamics procedure for finding the macroscale solution~$U_{j}(t_m)$ at time step~$t_m$, as well as subsequent time steps, is
\begin{enumerate}
\item Numerically solve equation~\eqref{eq:depatch} for the microscale fields within all patches~$j$ and for all configurations~$e$ at time~$t_m$ with coupling conditions~\eqref{eq:cc} where $U_{j}=U_{j}(t_{m-1})$ and $u_{j,i,e}=u_{j,i,e}(t_{m})$\,;\label{step:ujie}
\item From the amplitude condition~\eqref{eq:amp} determine the macroscale fields at time~$t_m$ at all macroscale lattice points~$X_j$ where both $u_{j,i,e}$ and~$U_{j}$ are at time~$t_m$; \label{step:U}
\item Using the macroscale solution~$U_{j}(t_m)$ evaluated at Step~\ref{step:U}, return to Step~\ref{step:ujie} to determine the solution at time step~$t_{m+1}$.
\end{enumerate} 
This process makes no attempt to derive an explicit macroscale model.
It simply evaluates `on the fly' numerical solutions of the model at macroscale lattice points~$X_j$ at each time step~$t_m$.
Here we only discuss macroscale modelling on a spatial scale and thus assume that standard numerical integrators are sufficiently accurate for evaluating $u_{j,i,e}(t_m)$ in Step~\ref{step:ujie} using previously calculated solutions at time step $t_{m-1}$.

\section{Optimising the patch geometry}
\label{sec:optimise}

We test the patch dynamics macroscale modelling by comparing the solution of equation~\eqref{eq:depatch} with coupling conditions~\eqref{eq:cc} to solutions of the original microscale discrete diffusion equation~\eqref{eq:de} on the complete domain.
One possibility is to compare field solutions $u_i(t)$ of the original model~\eqref{eq:de} at macroscale lattice points $X_j$ with patch dynamics field solutions $U_j(t)$ for various times~$t$.
However, we are primarily interested in how well patch dynamics captures the emergent large scale evolution of the macroscale field, as opposed to how well patch dynamics performs at fixed moments in time. 
Therefore, to determine the accuracy of patch dynamics we compare the large scale evolution~$g_0$, derived in Section~\ref{sec:infinite} on the complete domain, with the large scale evolution~$g_j$ derived from patch dynamics on the $j$th patch.

As discussed in Section~\ref{sec:infinite}, the large scale emergent evolution is characterised by the smallest magnitude eigenvalue~$\lambda_0$. 
The typical time for diffusion across the macroscale~$H$ is $H^2/4\kappa$, and on times scale of and beyond this diffusion time, the emergent dynamics should evolve as $g_0=\dot{U}_j=\lambda_0U_j/h^2$, as shown in equation~\eqref{eq:lambda0} derived on the complete domain. 
Similarly, we define patch dynamics large scale evolution as $g_j=\dot{U}_j$, but in this case $U_j$ is the macroscale field determined from patch dynamics.

Over the $j$th patch we solve for both the evolution~$g_j=g_j(\vec U)$ and the microscale fields~$u_{j,i,e}=u_{j,i,e}(\vec U)$ using the difference equation~\eqref{eq:depatch}, coupling conditions~\eqref{eq:cc} and amplitude condition~\eqref{eq:amp} with the time derivatives of the microscale fields defined by $\dot{u}_{j,i,e}(\vec U)=\sum_l g_l(\vec U)(du_{j,i,e}/dU_l)$ for macroscale index~$l$.
For a given~$j$ we need to solve for one evolution~$g_j$ and $2K(2n+1)$~microscale fields~$u_{j,i,e}$\,,  from $2K$ configurations~$e$ and $i=0,\pm 1,\ldots,\pm n$\,. We have a total of $2K(2n-1)$~difference equations,  two coupling conditions for each of the $2K$~configurations, and one amplitude condition. 
Analytic expressions are obtained for $g_j$ and $u_{j,i,e}$ in terms of the macroscale fields, the diffusivities $\kappa_i$, the diffusivity period~$K$, the patch half-width $n$ and the buffer half-width $b$.
The computer algebra package Reduce-Algebra\footnote{\url{http://reduce-algebra.com/}} efficiently solves this system of equations, although the computational time increases dramatically as $K$ or $n$ increase. The full code is freely available for download and modification~\cite{codevariD}.

Our aim is to show that patch dynamics is effective in determining the emergent evolution.
We determine the optimal patch geometry, defined by patch half-width $n$, and coupling conditions, defined by buffer half-width $b$, for a given period~$K$; that is, which combination of patch sizes and buffer sizes, $n$~and~$b$, produce an macroscale evolution~$g_j$ which is closest to the true macroscale evolution~$g_0$ derived in Section~\ref{sec:infinite} on the complete domain.
For a given diffusivity period~$K$, we evaluate $g_j$ for many patch half-widths~$n$ and buffer half-widths~$b$.
Section~\ref{sec:exact} discusses cases where certain patch geometries and coupling conditions produce an emergent macroscale evolution~$g_j$ which is precisely~$g_0$, to a specified order of accuracy.
These special cases are characterised by simple relationships between diffusivity period, patch geometry and coupling conditions, specifically  $K|(n-b)$~and~$K|(2b+1)$, and arise because symmetries of the original complete domain are not destroyed by the imposition of patches, unlike arbitrary $n$ and~$b$.
Thus, for a known~$K$ we should always aim to choose $n$ and~$b$ such that $g_j=g_0$ to a specified order of error.
However, there are circumstances where this may not be possible; 
for example, when the microscale period~$K$ is larger than our maximum choice of patch size~$n$, or when we are forced to use one particular~$n$.
Section~\ref{sec:notexact} considers the optimal buffer half-width~$b$ for a given patch half-width~$n$ for arbitrary diffusivity period~$K$ and finds that the best choice is generally a buffer width which is approximately the patch half-width, $b\approx n/2$\,.

Unless otherwise stated, to reduce the complexity of the computer algebra we assume there is only moderate variation between the diffusivities.
That is, we set $\kappa_i=\kappa_0(1+\eta_i)$ for some diffusivity $\kappa_0$ and dimensionless variation~$\eta_i$\,: the variations~\(\eta_i\) are big enough that quadratic effects are significant, but small enough that cubic effects are ignored. 
This simplification is not a requirement of our patch dynamics modelling but is implemented to avoid the need to deal with computationally expensive terms, such as rational terms with high order polynomials of $\vec\kappa$ in both the numerator and denominator.
In terms of small variations~$\vec\eta=(\eta_1,\ldots,\eta_K)$, the true macroscale, complete domain, evolution~\eqref{eq:lambda0} is 
\begin{align}
\dot U_j=g_0={}&\frac{\kappa_0 }{h^2}\left(1+d\sum_{i=1}^K\eta_i
-d_{0}\sum_{i=1}^K\eta_i^2+\sum_{i=1}^K\sum_{k=1}^{K/2}d_{k}
\eta_i\eta_{i+k}\right)\delta^2U_j\nonumber\\
{}&+\frac{\kappa_0 }{h^2}\left(f_0\sum_{i=1}^K\eta_i^2+\sum_{i=1}^K\sum_{k=1}^{K/2}f_{k}
\eta_i\eta_{i+k}\right)\delta^4U_j
+\mathcal{O}(\delta^6,\vec\eta^{\,3})\,, \label{eq:ex}
\end{align}
where the linear coefficient $d=1/K$\,, and the quadratic coefficients are
\begin{align}
d_{0}&=\frac{K-1}{K^2}\,,\quad 
d_{k}=\frac{2}{K}\text{ for }k\neq 0,K/2\,,\quad d_{K/2}=\frac{1}{K}\,,\quad f_0=\frac{K^2-1}{12K^2}\,,\nonumber\\
f_k&=\frac{1}{K^2}\left[\frac{K^2-1}{6}-k(K-k)\right]\text{ for }k\neq 0,K/2\,,\quad
f_{K/2}=-\frac{K^2+2}{24K^2}\,.
\end{align}
We find that patch dynamics reproduced the constant term and the linear term in~$\vec\eta$ exactly for a given~$K$ for all tested values of $n$~and~$b$. For our scheme, this remains true even when $K>n$ because the ensemble average ensures all diffusivites are accounted for, even when they cannot appear in one patch. 

To find the best performing patch coupling we explore effects quadratic in the diffusivity variations $\eta_i$ and the difference operator $\delta$, describe by the $d_k$ coefficients.
Such quadratic terms include the effects of physical correlations in the microscale structure that may be significant in applications (the \(\eta_i\eta_{i+k}\)~terms in~\eqref{eq:ex}). 
To compare the true macroscale evolution~$g_0$ with the patch dynamics macroscale evolution~$g_j$ we compare relative errors in the quadratic coefficients~$d_{k}$, 
\begin{equation}
\Delta d_{k}=\frac{d_k(\text{from }g_0)-d_k(\text{from }g_j)}{d_k(\text{from }g_0)}\quad \text{for }k=0,1,\ldots,K/2\,,
\label{eq:errdk}
\end{equation}
and we aim to determine what patch geometry and coupling conditions minimise the relative error~$\Delta d_k$\,.
We generally only investigate the~$\delta^2U_j$ term in the evolution expansion~\eqref{eq:ex} because this term is crucial to structurally stable macroscale models.
However, Section~\ref{sec:exact} discusses some special cases where the two macroscale descriptions, $g_0$~and~$g_j$, are identical, to the specified accuracy of~$g_0$, and so here we consider higher order coefficients $f_k$. 
In Section~\ref{sec:notexact} we consider the general case where $\Delta d_k$ is typically nonzero and, over the diffusivity period range $2\leq K\leq12$ with patch half-width range $2\leq n\leq 12$ and buffer half-widths $0\leq b<n$, we determine how to choose $n$ and $b$ such that $\Delta d_k$, and thus errors in $g_j$, are minimised.

%In these cases it is only realistic to expect patch dynamics agreement with the infinite domain solution up to first order in 
%the dimensionless diffusivity variability $\eta_i$\,,
%\begin{equation}
%g_0=\frac{\kappa}{h^2}\delta^2U_j+\mathcal{O}(\delta^4)=\frac{\kappa_a }{h^2}\delta^2U_j+\mathcal{O}(\delta^4,\eta^2)\,,
%\end{equation}  
%for $\kappa_a$ the arithmetic mean of the diffusivities in the sample of size~$K$. This is mainly a problem of
%choosing a appropriate sample size~$K$ for the given problem such that $\kappa_a$ 

\subsection{Ideal patch geometries give the correct macroscale}
\label{sec:exact}

When the microscale period~$K$ exactly divides the difference between the patch and buffer half-widths,~$K|(n-b)$, computer algebra~\cite{codevariD} shows that $g_j=g_0+\mathcal{O}(\delta^6,\vec\eta^{\,3})$ for all $2\leq K\leq12$ and $2\leq n\leq 12$\,.
Furthermore, for periods $K=2$ and~$3$ and~$K|(n-b)$, the equality $g_j=g_0+\mathcal{O}(\delta^6)$ holds for general diffusivities~$\kappa_i$\,; the accuracy of the patch dynamics modelling is not limited to cases with moderate variations in diffusivity $\kappa_i=\kappa_0(1+\eta_i)$ for small~$\eta_i$ and $g_0$ is defined by the more general equation~\ref{eq:lambda0} rather than equation~\ref{eq:ex}.
Computer algebra limitations prevent us from considering periods $K>3$ with general diffusivities. 

When the period~$K|(n-b)$, for the studied ranges of $K$ and $n$, to obtain $g_j=g_0+\mathcal{O}(\delta^6,\vec\eta^{\,3})$ we do not need the ensemble average.
To solve for the evolution~$g_j$ without an ensemble average we solve equation~\eqref{eq:depatch} for just one configuration~$e$ for $i=0,\pm 1,\ldots,\pm(n-1)$ on the $j$th~patch with coupling conditions~\eqref{eq:cc} where the macroscale fields~$U_j$ are derived from an amplitude condition like equation~\eqref{eq:amp}, but without the ensemble average.
In the single configuration case we only need to solve $2(n+1)$~equations, compared to $2K(2n+1)+1$~equations in the ensemble average case, leading to significant savings in the computation time required to solve the computer algebra, particularly if the period~$K$ is large. 

Reproducing $g_j=g_0+\mathcal{O}(\delta^6,\vec\eta^{\,3})$\,, with $g_0$ defined by equation~\ref{eq:ex}, is a significant achievement because it accurately describes the dependence of the macroscale evolution on the fine microscale detail.
Equation~\eqref{eq:ex} with error $\mathcal{O}(\delta^4,\vec\eta^{\,3})$, that is, ignoring all terms with $f_k$ coefficients, is expressible in terms of mean diffusivities (geometric, harmonic and algebraic) and is therefore not dependent on the fine microscale detail.
However,  with error $\mathcal{O}(\delta^6,\vec\eta^{\,3})$ equation~\ref{eq:ex} is no longer expressible as a function of mean diffusivities and is dependent on the detailed microscale structure, such as the ordering of the diffusivities on the microscale grid.
Thus, in accurately evaluating the coefficients $f_k$ patch dynamics provides an accurate macroscale closure of the microscale detail. 
The role of coefficients~$f_k$ in revealing the microscale structure of equation~\ref{eq:ex} is related to the role of $c_2(\vec \kappa)$ in equation~\eqref{eq:lambda0}, discussed in Section~\ref{sec:infinite}.
Both~$f_k$ and $c_2(\vec \kappa)$ show that the ordering of the diffusivities on the microscale lattice affects the macroscale evolution and a macroscale closure is not solely dependent on mean microscale quantities.

The previous paragraphs addressed modelling to errors~$\mathcal{O}(\delta^6)$.
We suspect that patch dynamics is accurate for errors of higher order than~$\mathcal{O}(\delta^6)$ since, for $K=2$ where the quadratic cutoff in Section~\ref{sec:infinite} is exact, higher order expansions of the square root in equation~\eqref{eq:lambda0} and computer algebra evaluations to the same order reveal $g_0$ and~$g_j$ are identical to any specified order of accuracy in arbitrary specified powers of~$\delta^2$. 

\begin{figure}
\centering\includegraphics{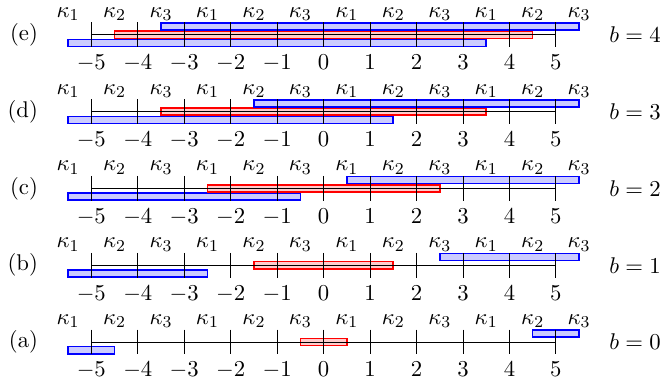}
\caption{Microscale structure of a single patch with half-width $n=5$\,, diffusivity period $K=3$ and different buffer half-widths $b$.
For each case, the core and buffers are indicated by the shaded outlined regions.
When $b=2$\,, $K|(n-b)$ is true and only then do the two buffers and core capture the same diffusivity pattern: $\kappa_1$~to the left, $\kappa_3$~to the right and 
$\kappa_2,\kappa_3,\kappa_1,\kappa_2$ inside.
When $b=(K-1)/2=1$ both buffer and core capture one complete diffusivity period, but not in the same order.
Similarly, when $b=(3K-1)/2=4$ both buffers and core capture two complete diffusivity periods.}
\label{fig:buffsymmetry}
\end{figure}

Patch dynamics is able to accurately predict the long scale evolution on the complete domain when~$K|(n-b)$ because the patch is reflecting the inherent symmetry of the original microscale model~\eqref{eq:de}. 
Our construction of the coupling conditions~\eqref{eq:cc} is reliant on a shift of $\pm(n-b)$ microscale lattice points of the microscale fields~$u_{j,i,e}$ in the patch core where $i=-b,\ldots,0,\ldots,b$ to the microscale fields in the buffers where $i=-b\pm(n-b),\ldots,\pm(n-b),\ldots,b\pm(n-b)$, as illustrated in Figure~\ref{fig:buffers}.
Therefore, since $K|(n-b)$, the core and both buffers contain identical diffusivities, as illustrated in Figure~\ref{fig:buffsymmetry}(c) for the case $K=3$ and $n=5$\,. By `identical' we mean that the three regions have identical diffusivities to the left, right and within.

As shown in Figure~\ref{fig:buffsymmetry}(b) and Figure~\ref{fig:buffsymmetry}(e), there are cases where the core and the buffers capture a multiple of a complete period of microscale structure, with one period captured in the former and two periods in the latter.
In general this requires $K|(2b+1)$ so is only possible when $K$~is odd. 
For~$K|(2b+1)$ we find that the patch dynamics evolution~$g_j$ accurately predicts the $d_{k}$~coefficients of the true macroscale evolution~$g_0$~\eqref{eq:ex} and $\Delta d_{k}=0$ for all~$k$.
However, this is only generally true when we include an ensemble average and moderate variation in diffusivities, which are not necessary when $K|(n-b)$\,. 
Also, unlike $K|(n-b)$, the case $K|(2b+1)$ does not produce an accurate evolution up to an arbitrary order of accuracy but produces errors in the quadratic coefficients~$f_k$ for all~$k$.
For even~$K$, and when using an ensemble average of moderately varying diffusivities,
although we cannot fit a multiple of diffusivity periods into the core or buffers, we do tend to find that
$\Delta d_k$ is relatively small for all~$k$ when $K|(2b+1)$ is approximately true, that is, $K|2b$ or~$K|2(b+1)$.

In general, when implementing patch dynamics, if we know the period~$K$ and are free to choose any patch geometry, then we should choose $n$~and~$b$ such that $K|(n-b)$ since it accurately predicts known emergent large scale evolution.
However, this accurate choice requires $K\leq n$\,, and when the microscale period~$K$ is large then the scheme may be computationally expensive.
For large~$K$, choosing the patch geometry such that $K|(2b+1)$ might appear to be a better choice since it is valid for smaller patch sizes $(K-1)/2<n$, but this case requires an ensemble average, which also becomes computationally expensive as~$K$ becomes large, albeit parallelisable.
Therefore, even when microscale period~$K$ is large we generally recommend choosing the patch geometry such that $K|(n-b)$ and evaluating over just one configuration, rather than an ensemble average.

\subsection{Optimum patch geometry when uncertain}
\label{sec:notexact}
Often we will not know the microscale period~\(K\).
This section explores the patch dynamics evolution~$g_j$ when the ideal geometry and coupling conditions satisfying~$K|(n-b)$ are not knowable. 
We seek a patch geometry and coupling conditions which are independent of the period~$K$ and minimise the error in second order coefficients in the patch dynamics evolution~$g_j$ of the macroscale evolution~\eqref{eq:ex}.
The error~\eqref{eq:errdk} is measured relative to the true macroscale evolution~$g_0$, defined by~$\Delta d_k$ for $k=0,\ldots,K/2$\,. 
We find that the general behaviour of the error~$\Delta d_k$ is both quantitatively and qualitatively similar for all~$k$, so here only present the $k=0$ case which measures the accuracy of~\(d_0\). 

\begin{figure}
\includegraphics{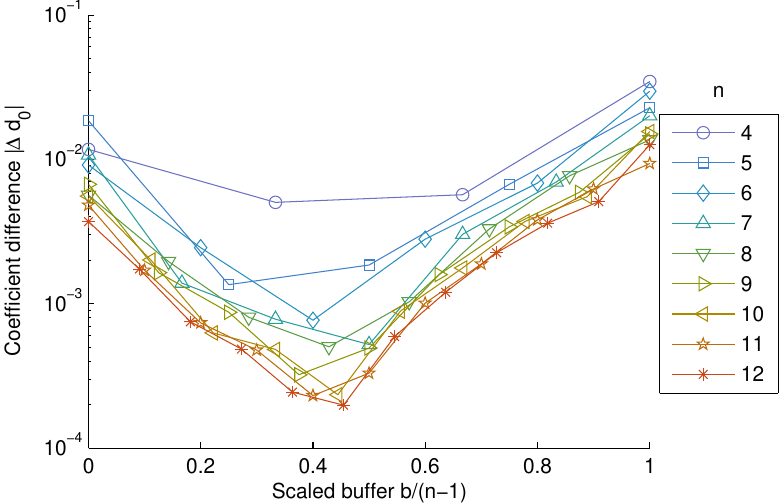}
\caption{Coefficient error~$|\Delta d_0|$ averaged over diffusivities $2\leq K\leq n$ relative to the scaled buffer $0\leq b/(n-1)\leq 1$ for $4\leq n\leq 12$\,. 
The coefficient error~$|\Delta d_0|$ is minimised when $b\approx n/2$\,, particularly when $n$ is large.}\label{fig:smallk0}
\end{figure}
 
Figure~\ref{fig:smallk0} plots errors~$|\Delta d_0|$ for given patch and buffers half-widths, $n$~and~$b$, averaged over all $2\leq K\leq n$\,. 
This figure caters for small microscale periods, those which fit within half a patch.
The worst choices are the smallest and largest possible buffers~$b=0,(n-1)$.
The error~$|\Delta d_0|$ decreases with increasing~$n$ and is minimised at $b\approx n/2$\,; that is, the optimal buffer size is generally half the patch size. 

\begin{figure}
\includegraphics{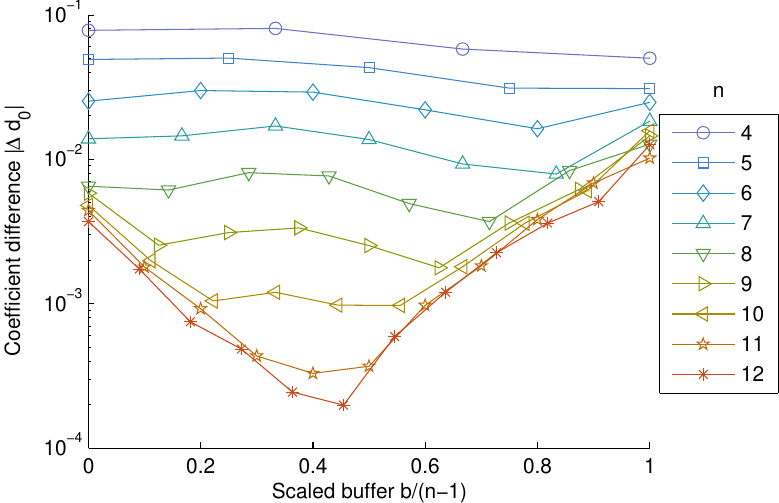}
\caption{Coefficient errors~$|\Delta d_0|$ averaged over diffusivities $2\leq K\leq 12$ relative to the scaled buffer $0\leq b/(n-1)\leq 1$ for $4\leq n\leq 12$\,. 
The coefficient error~$|\Delta d_0|$ is minimised when $b\approx n/2$ for large~$n$, but as $n$~decreases the minima are less obvious and shift to larger values of~$b/(n-1)$.}\label{fig:allk0}
\end{figure}

Figure~\ref{fig:allk0} plots the error~$|\Delta d_0|$ for given patch and buffers half-widths, $n$~and~$b$, but now averaged over a larger range of periods, namely $2\leq K\leq 12$\,. 
This average includes large periods which do not fit within one patch.
When most periods do not fit within one patch, as is the case for small~$n$, 
the difference~$|\Delta d_0|$ is rather large with slight minima near the largest buffer size $b\approx n-1$\,. 
As $n$~increases the minima become sharper and smaller as the patch is better able to accommodate all~$K$ diffusivities in the range $2\leq K\leq 12$\,. 
The position of the minima also shift to smaller buffer sizes and approaches the position $b\approx n/2$\,.

From Figure~\ref{fig:allk0} we conclude that patch dynamics does not accurately predict the macroscale evolution when the period of significant microscale structure,~$K$, 
is significantly larger than the patch half-width~$n$, with errors in the second order coefficients up to~$10\%$.
However, if $n$~is of the order of~$K$ or more, and the buffer half-width is chosen to be $b\approx n/2$\,, 
Figure~\ref{fig:smallk0} indicates that the evolution calculated from patch dynamics provides a good approximation to the evolution compared to the true macroscale evolution~$g_0$---the second order coefficients differing by less than~$1\%$.
For larger patches this error tends to decrease, with the errors in second order coefficients as small as~$0.1\%$ for the largest patches considered.

\section{Conclusion}
By analysing the emergent, large scale, evolution we showed that patch dynamics macroscale modelling is able to capture the emergent dynamics of a microscale lattice system when the microscale has fine detail.
For best results it is important to appropriately choose the patch geometry, defined by the patch half-width~$n$, and the patch coupling conditions, defined by the buffer half-width~$b$, relative to the underlying period~$K$ of the microscale detail.
We showed that the symmetry of the microscale model is important and should be reflected in the choice of patch geometry and coupling conditions.
We expect that other microscale models require similar consideration when implementing patch dynamics for macroscale modelling.

For the diffusion model considered here, the symmetry is best accounted for by choosing patch and buffer half-widths, $n$~and~$b$, such that microscale period~$K|(n-b)$ as not only is the complete domain evolution accurately obtained by patch dynamics, 
but the computational time can also be reduced by, without loss of accuracy,  
considering only one configuration rather than an ensemble average.
Furthermore, for this special case it appears the evolution is still accurate even when considering arbitrary microscale diffusivities. 
In cases where this ideal geometry is not available, best results are generally obtained when the patch half-width~$n$ is as large as possible and the buffer half-width is $b\approx n/2$ and when the macroscale fields are derived from an ensemble average. 

We defined~$K$ as the period of the microscale diffusivity, or equivalently the number of diffusivity values, but there are cases where a redefinition is required.
For example, when using raw data the microscale period may vary across all space and in this case~$K$ is simply redefined as some mean period and the patch geometry is chosen relative to this mean. 
A more significant example is when the period of the diffusivity is too large to be practical, or even infinite, as is effectively the case for a random medium.
In this case we redefine~$K$ as a sample size of diffusivity values and choose patch geometry relative to this~$K$.
How well patch dynamics is able to reproduce the long scale dynamics is dependent on the statistics of the sample size~$K$.

The microscale dynamics considered is left-right symmetric, since we only describe diffusion and not advection in the original microscale model.
We used this symmetry property to construct an ensemble of $2K$~configurations, $K$~translations and $K$~reflections of the original diffusivity configuration.
Future work will consider patch dynamics for macroscale modelling of microscale advection-diffusion equations.
The presence of advection terms will require a reconsideration of the ensemble.
%Certainly we cannot construct an ensemble which includes the~$K$ reflected configurations as these will cancel out the advection contributions.
%A suitable ensemble should account for the symmetries in the original microscale problem, and not introduce symmetries which do not exist in the original problem.
We expect that the patch geometry which best describes the emergent macroscale dynamics will complement any underlying symmetry of the detailed microscale problem.

Our analysis on the one dimensional, microscale variable, lattice diffusion equation~\eqref{eq:de} is readily modifiable to other systems.
For example, if the temporal derivative in equation~\eqref{eq:de} $\dot{u}_i$~is replaced with a second order derivative~$\ddot{u}_i$ we have a one dimensional lattice wave equation where~$\sqrt{\kappa_i}$ is a $K$~periodic, microscale varying, wave velocity.
The long scale dynamics is still described by the smallest magnitude eigenvalue~$\lambda_0$ but with $u_i\sim e^{\sqrt{\lambda_0}t/h}\sim U_j$ so that, over long scales, $\ddot{U}_j=\lambda_0/h^2$, which is comparable to~$\dot{U}_j$ for the diffusion equation, thus resulting in a very similar macroscale modelling problem.

\paragraph{Acknowledgements} The Australian Research Council Discovery Project grants DP0988738 and DP120104260 helped support this research.

%\bibliographystyle{plain} 
%\bibliography{../bibliography/bibliography}

\end{document}